\documentclass[12 pt]{article}
\usepackage{amsfonts}
\usepackage{amsmath}
\usepackage{pstricks}
\usepackage{graphicx}

\title{Sierpi\'nski Gasket Graphs and Some of Their Properties}
\author{Alberto M.~Teguia\\
Department of Mathematics\\
Duke University \and Anant P.~Godbole\\
Department of Mathematics\\
East Tennessee State University}
\begin{document}
\def\qed{\vbox{\hrule\hbox{\vrule\kern3pt\vbox{\kern6pt}\kern3pt\vrule}\hrule}}
\def\ep{\varepsilon}
\def\lr{\left(}
\def\lf{\lfloor}
\def\rf{\rfloor}
\def\lc{\left\{}
\def\rc{\right\}}
\def\dist{{\rm dist}}
\def\diam{{\rm diam}}
\def\rr{\right)}
\def\p{\mathbb P}
\def\e{\mathbb E}
\def\l{\lambda}
\def\lg{{\rm lg}}
\newtheorem{thm}{Theorem}
\newtheorem{lemma}[thm]{Lemma}
\newtheorem{prop}[thm]{Proposition}
\maketitle
\begin{abstract}  The {\it Sierpi\'nski fractal} or {\it Sierpi\'nski gasket} $\Sigma$ is a familiar object studied by specialists in dynamical systems and probability.  In this paper, we consider a graph $S_n$ derived from the first $n$ iterations of the process that leads to $\Sigma$, and study some of its properties, including its cycle structure, domination number and pebbling number.  Various open questions are posed.
\end{abstract}

\medskip \noindent {\it Keywords: Sierpi\'nski gasket; Sierpi\'nski gasket graph;  chromatic number; Hamiltonicity; pancyclicity; domination number; pebbling number; cover pebbling number.}

\section{Introduction and Basic Properties}  The structure known as the {\it Sierpi\'nski fractal} or {\it Sierpi\'nski gasket} $\Sigma$ is a familiar object studied by schoolchildren and specialists in dynamical systems and probability alike.  It is known to be a self-similar object with fractal dimension $d=\log3/\log2\approx1.585$ (\cite{edgar}).  In this paper, however, we consider the finite structure obtained by iterating, a finite number of times, the process that defines $\Sigma$, leading to (i) the finite Sierpi\'nski gasket $\sigma_n$, and (ii) the associated {\it Sierpi\'nski gasket graph} $S_n$, defined as one with vertex set $V_n$ equal to the intersection points of the line segments in $\sigma_n$, and edge set $E_n$ consisting of the line segments connecting two vertices; see Figure 1 for a portrayal of $S_2, S_3$ and $S_4$.  It is immediate that $S_{n+1}$ consists of three attached copies of $S_n$ which we will refer to as the {\it top, bottom left} and {\it bottom right} components of $S_{n+1}$ -- and denote by $S_{n+1,T}, S_{n+1,L},$ and $S_{n+1,R}$ respectively (Figure 2).

%Insert Figures 1 and 2 here
%begin figure 1

\unitlength 1mm

\psset{xunit=1mm,yunit=1mm,runit=1mm}
\psset{linewidth=0.3,dotsep=1,hatchwidth=0.3,hatchsep=1.5,shadowsize=1}
\psset{dotsize=0.7 2.5,dotscale=1 1,fillcolor=black}
\begin{center}
\begin{picture}(142.19,65.62)(0,0)
\linethickness{0.3mm} \put(11.25,30){\line(1,0){50}}
\linethickness{0.3mm}
\multiput(36.25,30)(0.12,0.17){104}{\line(0,1){0.17}}
\linethickness{0.3mm} \put(36.25,30.31){\circle*{1.88}}

\rput(36.25,25.31){$S_{2}$}

\rput(116.25,25.31){$S_{3}$}

\linethickness{0.3mm} \put(128.44,47.5){\circle*{1.88}}
\linethickness{0.3mm} \put(116.25,30){\circle*{1.88}}
\linethickness{0.3mm}
\multiput(36.25,65)(0.12,-0.17){208}{\line(0,-1){0.17}}
\linethickness{0.3mm}
\multiput(11.25,30)(0.12,0.17){208}{\line(0,1){0.17}}
\linethickness{0.3mm}
\multiput(23.75,47.5)(0.12,-0.17){104}{\line(0,-1){0.17}}
\linethickness{0.3mm} \put(23.75,47.5){\line(1,0){25}}
\linethickness{0.3mm} \put(141.25,30){\circle*{1.88}}
\linethickness{0.3mm} \put(11.56,30.31){\circle*{1.88}}
\linethickness{0.3mm} \put(48.44,47.5){\circle*{1.88}}
\linethickness{0.3mm} \put(36.25,64.69){\circle*{1.88}}
\linethickness{0.3mm} \put(24.06,47.5){\circle*{1.88}}
\linethickness{0.3mm} \put(61.25,30){\circle*{1.88}}
\linethickness{0.3mm} \put(91.25,30){\line(1,0){50}}
\linethickness{0.3mm}
\multiput(116.25,30)(0.12,0.17){104}{\line(0,1){0.17}}
\linethickness{0.3mm}
\multiput(116.25,65)(0.12,-0.17){208}{\line(0,-1){0.17}}
\linethickness{0.3mm}
\multiput(91.25,30)(0.12,0.17){208}{\line(0,1){0.17}}
\linethickness{0.3mm}
\multiput(103.75,47.5)(0.12,-0.17){104}{\line(0,-1){0.17}}
\linethickness{0.3mm} \put(103.75,47.5){\line(1,0){25}}
\linethickness{0.3mm} \put(104.06,47.5){\circle*{1.88}}
\linethickness{0.3mm} \put(116.25,64.69){\circle*{1.88}}
\linethickness{0.3mm} \put(122.81,38.75){\circle*{1.88}}
\linethickness{0.3mm} \put(91.56,30){\circle*{1.88}}
\linethickness{0.3mm} \put(97.5,38.75){\line(1,0){12.5}}
\linethickness{0.3mm}
\multiput(97.5,38.75)(0.12,-0.17){52}{\line(0,-1){0.17}}
\linethickness{0.3mm}
\multiput(103.75,30)(0.12,0.17){52}{\line(0,1){0.17}}
\linethickness{0.3mm} \put(110,56.25){\line(1,0){12.5}}
\linethickness{0.3mm}
\multiput(110,56.25)(0.12,-0.17){52}{\line(0,-1){0.17}}
\linethickness{0.3mm}
\multiput(122.5,38.75)(0.12,-0.17){52}{\line(0,-1){0.17}}
\linethickness{0.3mm}
\multiput(116.25,47.5)(0.12,0.17){52}{\line(0,1){0.17}}
\linethickness{0.3mm}
\multiput(128.75,30)(0.12,0.17){52}{\line(0,1){0.17}}
\linethickness{0.3mm} \put(122.5,38.75){\line(1,0){12.5}}
\linethickness{0.3mm} \put(135,38.75){\circle*{1.88}}
\linethickness{0.3mm} \put(128.75,30.31){\circle*{1.88}}
\linethickness{0.3mm} \put(116.25,47.81){\circle*{1.88}}
\linethickness{0.3mm} \put(122.5,56.25){\circle*{1.88}}
\linethickness{0.3mm} \put(110.31,55.94){\circle*{1.88}}
\linethickness{0.3mm} \put(97.5,38.75){\circle*{1.88}}
\linethickness{0.3mm} \put(103.75,30){\circle*{1.88}}
\linethickness{0.3mm} \put(110,38.75){\circle*{1.88}}
\end{picture}

\begin{pspicture}(0,0)(152.5,67.5)
\rput{0}(98.75,5){\psellipse[fillstyle=solid](0,0)(0.94,0.94)}
\rput{0}(123.75,5){\psellipse[fillstyle=solid](0,0)(0.94,0.94)}
\psline(30,-3.75)(90,-3.75) \psline(92.5,48.75)(79.5,31)
\psline(105,31.25)(130,-3.75) \psline(55,31.25)(80.5,66.5)
\psline(42.5,13.75)(55,-3.75) \psline(92.5,48.75)(67.5,48.75)
\rput{0}(61.25,39.68){\psellipse[fillstyle=solid](0,0)(0.94,0.94)}
\rput{0}(86.25,57.5){\psellipse[fillstyle=solid](0,0)(0.94,0.94)}
\rput{0}(105,31.56){\psellipse[fillstyle=solid](0,0)(0.94,0.94)}
\rput{0}(42.5,-3.44){\psellipse[fillstyle=solid](0,0)(0.94,0.94)}
\psline[fillstyle=solid](86.25,40)(98.75,40)
\psline[fillstyle=solid](36.25,5)(42.5,-3.75)
\psline[fillstyle=solid](61.25,22.5)(55,13.75)
\psline[fillstyle=solid](73.75,57.5)(86.25,57.5)
\psline[fillstyle=solid](61.25,40)(67.5,31.25)
\psline[fillstyle=solid](48.75,22.5)(55,13.75)
\psline[fillstyle=solid](73.75,40)(67.5,31.25)
\psline[fillstyle=solid](48.75,5)(42.5,-3.75)
\psline[fillstyle=solid](48.75,22.5)(61.25,22.5)
\rput{0}(105,-3.44){\psellipse[fillstyle=solid](0,0)(0.94,0.94)}
\rput{0}(117.5,-3.75){\psellipse[fillstyle=solid](0,0)(0.94,0.94)}
\rput{0}(86.25,40){\psellipse[fillstyle=solid](0,0)(0.94,0.94)}
\rput{0}(67.5,48.75){\psellipse[fillstyle=solid](0,0)(0.94,0.94)}
\rput{0}(74.06,57.5){\psellipse[fillstyle=solid](0,0)(0.94,0.94)}
\rput{0}(67.5,31.56){\psellipse[fillstyle=solid](0,0)(0.94,0.94)}
\rput{0}(80,31.56){\psellipse[fillstyle=solid](0,0)(0.94,0.94)}
\rput{0}(105,14.06){\psellipse[fillstyle=solid](0,0)(0.94,0.94)}
\rput{0}(98.75,22.5){\psellipse[fillstyle=solid](0,0)(0.94,0.94)}
\rput{90}(129.68,-3.44){\psellipse[fillstyle=solid](0,0)(0.94,0.94)}
\psline(55,31.25)(105,31.25) \psline(67.5,13.75)(55,-3.75)
\psline(80,66.25)(105,31.25) \psline(30,-3.75)(55,31.25)
\psline(67.5,48.75)(80,31.25) \psline(67.5,13.75)(42.5,13.75)
\rput{0}(73.75,40){\psellipse[fillstyle=solid](0,0)(0.94,0.94)}
\rput{0}(80,65.94){\psellipse[fillstyle=solid](0,0)(0.94,0.94)}
\rput{0}(30,-3.44){\psellipse[fillstyle=solid](0,0)(0.94,0.94)}
\psline[fillstyle=solid](36.25,5)(48.75,5)
\psline[fillstyle=solid](86.25,40)(92.5,31.25)
\psline[fillstyle=solid](98.75,40)(92.5,31.25)
\psline[fillstyle=solid](61.25,40)(73.75,40)
\psline[fillstyle=solid](73.75,57.5)(80,48.75)
\psline[fillstyle=solid](61.25,5)(67.5,-3.75)
\psline[fillstyle=solid](86.25,57.5)(80,48.75)
\psline[fillstyle=solid](73.75,5)(67.5,-3.75)
\psline[fillstyle=solid](61.25,5)(73.75,5)
\rput{0}(117.5,13.75){\psellipse[fillstyle=solid](0,0)(0.94,0.94)}
\rput{0}(111.25,5){\psellipse[fillstyle=solid](0,0)(0.94,0.94)}
\rput{0}(98.44,40.32){\psellipse[fillstyle=solid](0,0)(0.94,0.93)}
\rput{0}(92.5,48.75){\psellipse[fillstyle=solid](0,0)(0.94,0.94)}
\rput{0}(80,49.06){\psellipse[fillstyle=solid](0,0)(0.94,0.94)}
\rput{90}(55.32,31.56){\psellipse[fillstyle=solid](0,0)(0.94,0.94)}
\rput{90}(92.18,13.75){\psellipse[fillstyle=solid](0,0)(0.94,0.94)}
\rput{0}(92.5,31.56){\psellipse[fillstyle=solid](0,0)(0.94,0.94)}
\psline(117.5,13.75)(105,-3.75) \psline(55,31.25)(80,-3.75)
\psline(117.5,13.75)(92.5,13.75)
\rput{0}(61.56,5){\psellipse[fillstyle=solid](0,0)(0.94,0.94)}
\rput{0}(55,14.06){\psellipse[fillstyle=solid](0,0)(0.94,0.94)}
\rput{0}(111.25,22.5){\psellipse[fillstyle=solid](0,0)(0.94,0.94)}
\psline[fillstyle=solid](86.25,5)(98.75,5)
\psline[fillstyle=solid](98.75,22.5)(111.25,22.5)
\psline[fillstyle=solid](98.75,22.5)(105,13.75)
\psline[fillstyle=solid](98.75,5)(92.5,-3.75)
\rput{0}(92.5,-3.44){\psellipse[fillstyle=solid](0,0)(0.94,0.94)}
\rput{0}(42.5,13.44){\psellipse[fillstyle=solid](0,0)(0.94,0.94)}
\rput{90}(49.06,22.5){\psellipse[fillstyle=solid](0,0)(0.94,0.94)}
\rput{0}(36.56,5){\psellipse[fillstyle=solid](0,0)(0.94,0.94)}
\rput{0}(55,-3.44){\psellipse[fillstyle=solid](0,0)(0.94,0.94)}
\psline(80,-3.75)(130,-3.75) \psline(80,-3.75)(105,31.25)
\psline(92.5,13.75)(105,-3.75)
\rput{0}(73.44,5){\psellipse[fillstyle=solid](0,0)(0.94,0.94)}

\rput{0}(60.94,22.5){\psellipse[fillstyle=solid](0,0)(0.94,0.94)}
\psline[fillstyle=solid](111.25,5)(117.5,-3.75)
\psline[fillstyle=solid](123.75,5)(117.5,-3.75)
\psline[fillstyle=solid](111.25,5)(123.75,5)
\psline[fillstyle=solid](86.25,5)(92.5,-3.75)
\psline[fillstyle=solid](111.25,22.5)(105,13.75)
\rput{90}(80,-3.44){\psellipse[fillstyle=solid](0,0)(0.94,0.94)}
\rput{0}(86.56,5){\psellipse[fillstyle=solid](0,0)(0.94,0.94)}
\rput{0}(67.5,13.43){\psellipse[fillstyle=solid](0,0)(0.94,0.94)}
\rput{0}(48.44,5){\psellipse[fillstyle=solid](0,0)(0.94,0.94)}
\rput{0}(67.5,-3.44){\psellipse[fillstyle=solid](0,0)(0.94,0.94)}
\rput(50,22.5){} \rput(77.5,67.5){1} \rput(90,57.5){3}
\rput(65,50){4} \rput(80,47.5){5} \rput(95,50){6} \rput(57.5,40){7}
\rput(77.5,40){8} \rput(83.5,40){9}\rput(102.5,40){10}
\rput(70,57.5){2} \rput(52.5,32.5){11} \rput(67.5,27.5){12}
\rput(80,27.5){13} \rput(92.5,27.5){14} \rput(107.5,32.5){15}
\rput(45,22.5){16} \rput(65,22.5){17} \rput(95,22.5){18}
\rput(115,22.5){19} \rput(37.5,15){20} \rput(55,10){21}
\rput(87.5,15){23} \rput(72.5,15){22} \rput(105,10){24}
\rput(122.5,15){25} \rput(32.5,5){26} \rput(52.5,5){27}
\rput(57.5,5){28} \rput(77.5,5){29} \rput(82.5,5){30}
\rput(102.5,5){31} \rput(107.5,5){32} \rput(127.5,5){33}
\rput(25,-2.5){34} \rput(42.5,0){35} \rput(55,0){36}
\rput(67.5,0){37} \rput(80,0){38} \rput(92.5,0){39} \rput(105,0){40}
\rput(117.5,0){41} \rput(132.5,-2.5){42} \rput(80,-10){$S_{4}$}
\rput(80,-30){Figure 1}
\end{pspicture}
\end{center}
%end figure 1

%begin figure2

\psset{xunit=1mm,yunit=1mm,runit=1mm}
\psset{linewidth=0.3,dotsep=1,hatchwidth=0.3,hatchsep=1.5,shadowsize=1}
\psset{dotsize=0.7 2.5,dotscale=1 1,fillcolor=black}
\begin{flushleft}
\begin{pspicture}(0,0)(141.16,114)
\psline(20,30)(140,30) \psline(20,30)(80,110)
\psline(140,30)(80,110) \psline(50,70)(80,30) \psline(80,30)(110,70)
\psline(110,70)(50,70) \rput(80,80){$S_{n+1,T}$}
\rput(110,40){$S_{n+1,R}$} \rput(50,41){$S_{n+1,L}$}
\rput(80,-2){Figure 2, $S_{n+1}$} \rput(133,74){}
\rput(133,74){$S_{n+1,T,R}=S_{n+1,R,T}$} \rput(80,20){}
\rput(80,20){$S_{n+1,L,R}=S_{n+1,R,L}$} \rput(20,24){$S_{n+1,L,L}$}
\rput(22,72){} \rput(21,72){$S_{n+1,T,L}=S_{n+1,L,T}$}
\rput(80,114){$S_{n+1,T,T}$} \rput(140,25){$S_{n+1,R,R}$}
\rput{90}(50,70){\psellipse[fillstyle=solid](0,0)(1.58,1.58)}
\rput{90}(110,70){\psellipse[fillstyle=solid](0,0)(1.58,1.58)}
\rput{90}(139.58,30){\psellipse[fillstyle=solid](0,0)(1.58,1.58)}
\rput{90}(80,30){\psellipse[fillstyle=solid](0,0)(1.58,1.58)}
\rput{90}(20.42,30){\psellipse[fillstyle=solid](0,0)(1.58,1.58)}
\rput{90}(80,109.58){\psellipse[fillstyle=solid](0,0)(1.58,1.58)}
\end{pspicture}
\end{flushleft}
%end figure 2

It turns out that structures similar to ours have been studied in two other contexts.  It is unfortunate, but not surprising, that the phrase ``Sierpi\'nski graph" is used in each of these two situations, described below, as well.  There will be little danger of confusion, however, since we will deal exclusively, in this paper, with the definition in the previous paragraph.

First, we mention the body of work that treats simple random walks and Brownian motion on the ``infinite Sierpi\'nski graph" and the Sierpi\'nski gasket respectively.  A typical publication in this genre is by Teufl \cite{teufl}, where a sharp average displacement result is proved for random walk on a Sierpi\'nski graph that is an infinite version of the graph $S_n$  defined above -- but with each edge being of length one and with each vertex having degree four -- in our case all but three vertices have degree four and the edges can be arbitrarily short.  Several key references may be found in Teufl's paper.

Second, in \cite{mohar}, for example, the authors study crossing numbers for ``Sierp\-i\'n\-ski-like graphs," which are extended versions of the ``Sierpi\'nski graphs" $S(n,k)$ first studied in \cite{klavzar}.   These graphs were motivated by topological studies of the Lipscomb space that generalizes the Sierpi\'nski gasket and are defined as follows:  $S(n,k)$ has vertex set $\{1,2,\ldots,k\}^n$, and there is an edge between two vertices $u=(u_1,\ldots,u_n)$ and $v=(v_1,\ldots,v_n)$ iff there is an $h\in[n]$ such that
\begin{itemize}
\item $u_j=v_j$ for $j=1,\ldots,h-1$;
\item $u_h\ne v_h$; and
\item $u_j=v_h; v_j=u_h$ for $j=h+1,\ldots,n$.
\end{itemize}
An interesting connection is that the graph $S(n,3),n\ge1$ is isomorphic to the Tower of Hanoi game graph with $n$ disks:  see \cite{touhey}
for a popular account and \cite{klavzarhanoi} 
for a list of scholarly references.

After discussing some baseline results in this section, we will prove in Section 2 that $S_n$ is Hamiltonian and pancyclic (i.e., has cycles of all possible sizes).  In Section 3, we determine the domination number of $S_n$, proving in the process that it has efficiency that is asymptotically 90\%.  This is in sharp contrast to the Sierpi\'nski graphs of \cite{klavzar} -- it has been exhibited by Klav\v zar, Milutinovi\'c, and Petr \cite{klavzar2} that the graphs $S(n,k)$ have perfect dominating sets, i.e., are 100\% efficient.  Finally, in Section 4, we will show that the determination of the cover pebbling number $\pi(S_n)$ of $S_n$ is not trivial, even given the so-called ``stacking theorem" of Vuong and Wyckoff \cite{ian}.  We conclude with a statement of some open questions.
\begin{prop}
$S_n$ has ${3\over2}\lr3^{n-1}+1\rr$ vertices and $3^n$ edges. 
\end{prop}

\medskip

\noindent{\bf Proof} To
construct $S_{n+1}$ from $S_n$ we add one downward facing triangle in
each of the $3^{n-1}$ upward facing triangles of $ S_n $. Thus we
add $3^n$ points.  In other words,

\begin{eqnarray}
  |V_{n+1}| &=& |V_n| + 3^n \nonumber \\
  &=& |V_1| + \sum_{i=1}^n 3^i\nonumber\\
  &=&\frac{3}{2}\cdot(3^n + 1),
  \nonumber
\end{eqnarray} as asserted.  The number of edges in $S_n$ may now be easily  determined using the fact that the sum of the vertex degrees equals twice the number of edges:
\begin{eqnarray*}
\vert E_n\vert&=&\frac{1}{2}\sum_{j=1}^{\vert V_n\vert}\deg(v_j)\\
&=&\frac{1}{2}\cdot\lr4\cdot\frac{3}{2}(3^{n-1}-1)+2\cdot3\rr\\
&=&3^n,
\end{eqnarray*}
completing the proof.\hfill\qed

Note that an alternative proof of Proposition 1 can be based on the facts that $\vert V_{n+1}\vert=3\vert V_n\vert-3$ and $\vert E_{n+1}\vert=3\vert E_n\vert.$
  
\begin{prop} $S_n$ is properly three-colorable, i.e. $\chi(S_n)=3$ for each $n$.
\end{prop}
\medskip

\noindent{\bf Proof}  Clearly, $\chi(S_1) = 3$.  Suppose $\chi(S_n)= 3$.  Color $S_n$ with three colors.  We now properly 3-color $S_{n+1}$ in the following fashion:  After we insert one downward facing triangle in each upward facing
triangle of $S_n$, each added vertex
is assigned a previously used color different from
the colors of the vertices of $S_n$ adjacent to it.
\hfill\qed

\section{Hamiltonicity and Pancyclicity}  We begin with an important lemma:
\begin{lemma}
$S_n$ has two Hamiltonian paths, say $H_{n0}$ and $H_{n1}$, both
starting
 at the same vertex of degree two and ending at different vertices of degree two.
\end{lemma}

\medskip

\noindent{\bf Proof} The lemma is clearly true for $S_1$.
Suppose it is true for $S_n$.  Consider $S_{n+1}$, which is a three-fold ``repetition" of $S_n$, as mentioned in Section 1, i.e., it consists of the three attached copies $S_{n+1,T}, S_{n+1,L},$ and $S_{n+1,R}$ of $S_n$.  Consider a  Hamiltonian path of $S_{n+1,T}$, moving from the
top vertex (of degree 2) of $S_{n+1,T}$ to the top vertex of $S_{n+1,R}$.  Using another Hamiltonian
path (guaranteed by the induction hypothesis) we move from that vertex
to left vertex of of $S_{n+1,R}$.
Finally, we take the Hamiltonian path of $S_{n+1,L}$ starting at its right vertex and ending at its left vertex, but with a critical modification, namely avoiding the top vertex of $S_{n+1,L}$.  In this fashion, we have constructed a Hamiltonian path of
$S_{n+1}$ from its top vertex to its left vertex.  A similar argument is employed if the path is to end in the right vertex of $S_{n+1}$.
\hfill\qed

We simplify our notation next.  The top, left, and right vertices of $S_{n+1,T}$ will be denoted respectively by $S_{n+1,T,T},$ $S_{n+1,T,L},$ and $S_{n+1,T,R}$.  Other critical vertices in $S_{n+1}$ are analogously denoted by $S_{n+1,L,T},S_{n+1,L,L},S_{n+1,L,R},$

\noindent $S_{n+1,R,T},S_{n+1,R,L},$ and $S_{n+1,R,R}$; of course we have $S_{n+1,T,R}=S_{n+1,R,T};$

\noindent $S_{n+1,R,L}=S_{n+1,L,R};$ and $S_{n+1,L,T}=S_{n+1,T,L}$.  See Figure 2.

\begin{thm} $S_n$ is Hamiltonian for each $n$.
\end{thm}

\medskip

\noindent{\bf Proof}
By Lemma 3, we take a Hamiltonian path of $S_{n+1,T}$ that
moves from $S_{n+1,T,L}$ to $S_{n+1,T,R}$. Next, take a Hamiltonian path from $S_{n+1,T,R}$ to $S_{n+1,R,L}$ and finally move from this vertex to our
starting vertex using another Hamiltonian path of $S_{n}$.
This gives us a Hamiltonian cycle for $S_{n+1}$.
\hfill\qed
\begin{lemma}
Each  Hamiltonian path of $S_{n+1}$ as constructed in Lemma 3 above, and moving, say, from $S_{n+1,T,T}$ to $S_{n+1,L,L}$, can be sequentially
reduced in length by one at each step, while maintaining the starting and ending vertices, with the process ending in a path from $S_{n+1,T,T}$ to $S_{n+1,L,L}$ along the left side of  $S_{n+1}$.
\end{lemma}

\medskip

\noindent{\bf Proof}  We proceed by induction, noting that the result is clearly true for $n=1$.  Assume that the result is true for $S_n$.  There exists a Hamiltonian path from $S_{n+1,T,T}$ to $S_{n+1,L,L}$ via $S_{n+1,T,R}$ and $S_{n+1,R,L}$.  Suppose we need a path from $S_{n+1,T,T}$ to $S_{n+1,L,L}$ but with a reduction in length of $r\ge1$.  If $r\le\vert S_n\vert-s$, where $s$ is the length of the side of $S_n$, then adjustments need only be made in the first of the three above-mentioned Hamiltonian paths.  If $\vert S_n\vert-s\le r\le2\vert S_n\vert-2s$, then we make adjustments in the lengths of two Hamiltonian paths.  Likewise, if $2\vert S_n\vert-2s\le r\le3\vert S_n\vert-3s$, we adjust three paths.  Finally, if $3\vert S_n\vert-3s\le r\le3\vert S_n\vert-2s$, then we modify the reduced path $S_{n+1,T,T}\rightarrow S_{n+1,T,R}\rightarrow S_{n+1,R,L}\rightarrow S_{n+1,L,L}$ of length $3s$ as follows to achieve the further required reduction:  (i) $S_{n+1,T,T}\rightarrow S_{n+1,T,R}\rightarrow S_{n+1,R,L}$ is replaced by a path from $S_{n+1,T,T}$ to $S_{n+1,T,L}$ of length $2s$; (ii) the path $S_{n+1,R,L}\rightarrow S_{n+1,L,L}$ of length $s$ is replaced by the path $S_{n+1,T,L}\rightarrow S_{n+1,L,L}$, also of length $s$.  This yields a path from $S_{n+1,T,T}$ to $S_{n+1,L,L}$ of length $3s$. Finally, the first component of this path is shrunk, by the induction hypothesis, to a path of length $s+\delta$, leading to a path of length $2s+\delta; 0\le\delta\le s$ from $S_{n+1,T,T}$ to $S_{n+1,L,L}$.  This completes the proof.\hfill\qed
\begin{thm} $S_n$ is pancyclic for each $n$.
\end{thm}

\medskip

\noindent{\bf Proof}  Once again the proof is by induction.  Assume the result is true for all $m\le n$.  The Hamiltonian cycle $S_{n+1,T,L}\rightarrow S_{n+1,T,R}\rightarrow S_{n+1,R,L}\rightarrow S_{n+1,T,L}$ of $S_{n+1}$ consists of three Hamiltonian paths in $S_{n+1,T}$, $S_{n+1,R}$ and $S_{n+1,L}$ respectively.  By Lemma 5 applied to $S_n$, we reduce these as necessary to get cycles of all sizes $\ge3s$.  Cycles of smaller sizes are obtained by invoking the induction hypothesis on $S_{n}$, noting that $\vert S_n\vert\ge3s$.\hfill\qed

\section{Domination Numbers and Efficiency}  The three degree two vertices $\{v_T,v_L,v_R\}$ of $S_n$
will be called ``corner vertices" and the three vertices $S_{n,T,R},S_{n,R,L}$ and $S_{n,L,T}$ will be called the ``middle vertices" of $S_n$.
Let $S_n' = S_n\setminus\{v_T,v_L,v_R\}$. For $k=0,1,2,3$, let ${\gamma_n}^k$ be the
minimum number of vertices needed to dominate
$S_n'$ {\it in addition to} $k$ corner vertices that are to assist in the dominating of $S_n'$.  Let $\gamma_n$ denote the domination number of $S_n$.

\begin{thm}
 For every $n  \ge4$ we have
\begin{displaymath}
\gamma_n = 3 \cdot \gamma_{n-1}
\end{displaymath}

and

\begin{displaymath}
{\gamma_n}^k \ge\left\{ \begin{array}{ll}
\vspace{0.05in}
~\gamma_n & \textrm{if $k$ = 0, 1} \\
\vspace{0.05in}
\gamma_n -1 & \textrm{if $k$ = 2, 3} \\

\end{array} \right.
\end{displaymath}
\end{thm}

\medskip

\noindent{\bf Proof}  The fact that $\gamma_{n}\le3\gamma_{n-1}$ is a trivial consequence of the decomposition of $S_{n+1}$ into its three components.  

It is immediate that $\gamma_1=1;\gamma_2=2;$ and $\gamma_3=3.$ We next verify that both parts of the result are true for $n=4$.  First note that $\gamma_4\le9$.  Also if $\gamma_4=8,$ we contradict the fact that for any graph $G$ with maximum degree $\Delta$, $\vert\gamma(G)\vert(\Delta+1)\ge\vert V(G)\vert$.  The fact that $\gamma_4^0=\gamma_4^1=9$ may be checked by hand.  Consider $\gamma_4^2$.  Since a total of four vertices of $S_4'$ are dominated by the two external vertices, we need to dominate 35 others in $S_4'$.  Assuming that the two aiding vertices are $v_1$ and $v_{34}$ (see Figure 1), we must have $v_4$ or $v_6$ and $v_{20}$ or $v_{36}$ in the dominating set for $S_4'$.  Supposing without loss of generality that $v_4$ and $v_{20}$ are in the dominating set, we must now dominate 27 additional vertices and thus need at least 6 other vertices in the dominating set. Thus $\gamma_4^2\ge 8$ as required.  The fact that $\gamma_4^3\ge8$ is checked similarly.  

Assume then that the statements of the theorem are both true for each $m$ with $4\le m\le n$.
Let us start by proving the first part of the theorem.  Since any dominating set of $S_{n+1}$ contains either 0, 1, 2, or 3 ``middle vertices," we have
\begin{eqnarray}
  \gamma_{n+1} &\geq& \min \{3 {\gamma_n}^0, 2
{\gamma_n}^1+{\gamma_n}^0+1, 2
 {\gamma_n}^1 + {\gamma_n}^2+2, 3
{\gamma_n}^2+3 \}  \nonumber \\
  &\ge& \min \{3\gamma_n, 3
\gamma_n+1, 3\gamma_n+1, 3\gamma_n\} \nonumber\\
  &=& 3\gamma_n,
 \nonumber
\end{eqnarray}
as required.  A word of explanation might be in order:  In the above calculation, the first quantity, namely $3\gamma_n^0$, is a lower bound on $\gamma_{n+1}$ assuming that no middle vertices are in the dominating set of $S_{n+1}$.  It is obtained as follows.  There might be 0, 1, 2, or 3 {\it corner} vertices in the dominating set and we thus have, in this case, 
\[\gamma_{n+1}\ge\min\{3\gamma_n^0,2\gamma_n^0+\gamma_n^1+1, \gamma_n^0+2\gamma_n^1+2, 3\gamma_n^1+3\}=3\gamma_n^0.\]
Actually, it is evident that the minimum in each case corresponds to there being {\it no} corner vertices in the dominating set.

 For the second part of the proposition, we note that $\gamma_n^0\ge\gamma_n^1$ and $\gamma_n^2\ge\gamma_n^3$, so that

\begin{eqnarray}
 {\gamma_{n+1}^0\ge\gamma_{n+1}}^1 &\geq&\min \{2  {\gamma_n}^0+{\gamma_n}^1,
  3{\gamma_n}^1+1,{\gamma_n}^2+
{\gamma_n}^1+{\gamma_n}^0+1 ,\nonumber\\
& & {\gamma_n}^3 + 2
{\gamma_n}^1 + 2 , {\gamma_n}^1 + 2
 {\gamma_n}^2 + 2 , {\gamma_n}^3
+ 2
 {\gamma_n}^2 + 3  \}  \nonumber \\
&=&3\gamma_n\nonumber\\
  &=&\gamma_{n+1},
  \nonumber
\end{eqnarray}
and
\begin{eqnarray*}
 {\gamma_{n+1}^2\ge\gamma_{n+1}}^3 &\geq &\min \{3
 {\gamma_n}^1,2{\gamma_n}^2 +
 {\gamma_n}^1 + 1,2{\gamma_n}^2 +
 {\gamma_n}^3 + 2,3{\gamma_n}^3+3\}\\
 &\ge& 3 \gamma_n-1\\
&=&\gamma_{n+1}-1,
\end{eqnarray*}
completing the proof.   \hfill\qed

\medskip

\noindent{\bf Remarks}  Note that by Theorem 7, $\gamma_n=3^{n-2}$ for $n\ge3$, and thus, since none of the outer vertices are in the minimum dominating set, it follows that this set ``covers," with duplication, a total of $5\cdot3^{n-2}$ vertices.  Now $S_n$ has ${3\over2}\lr3^{n-1}+1\rr$ vertices, so that the ``efficiency" of the domination is asymptotically 90\%.  This is in contrast to the fact, exhibited by Klav\v zar, Milutinovi\'c, and Petr \cite{klavzar2} that the graphs $S(n,k)$ have perfect dominating sets, i.e., are 100\% efficient.  After we had completed this research, our colleague Teresa Haynes pointed out that domination numbers of the so-called E-graphs of \cite{teresa} provide generalizations of Theorem 7.  For this reason, we have given only an abbreviated proof of Theorem 7.

\section{Pebbling Numbers}
Given a connected graph $G$, distribute $t$ indistinguishable pebbles on its vertices in some
configuration. Specifically, a configuration of weight $t$ on a graph $G$ is
a function $C$ from the vertex set $V(G)$ to $\mathbf{N} \cup \{0\} $ such that $\sum_{v\in V(G)}C(v)=t$.
A {\it pebbling move} is defined as the removal of two pebbles from
some vertex and the placement of one of these on an adjacent vertex. Given an initial configuration, a vertex $v$ is called {\it reachable} if it is possible to place a pebble on it in finitely many pebbling moves.  Given a configuration, the graph $G$ is said to be {\it pebbleable} if any of its vertices can be thus reached. Define the pebbling
number $\pi(G)$ to be the minimum number of pebbles that are sufficient to pebble the graph regardless of the initial configuration.

\medskip

\noindent SPECIAL CASES:  The pebbling number $\pi(P_n)$ of the path is $2^{n-1}$ (\cite{glennsurvey}). Chung \cite{chung} proved that $\pi(Q^d)=2^d$ and $\pi(P_n^m)=2^{(n-1)m}$, where $Q^d$ is the $d$-dimensional binary cube and $P_n^m$ is the cartesian product of $m$ copies of $P_n$.    An easy pigeonhole principle argument yields $\pi(K_n)=n$.  The pebbling number of trees has been determined (see \cite{glennsurvey}).
One of the key conjectures in pebbling, now proved in several special cases, is due to Graham; its resolution would clearly generalize Chung's result for $m$-dimensional grids:

\medskip

\noindent GRAHAM'S CONJECTURE.  {\it The pebbling number of the cartesian product of two graphs is no more that the product of the pebbling numbers of the two graphs, i.e.
$$\pi(G \Box H)\le \pi(G)\pi(H).$$}

A detailed survey of graph pebbling has been presented by
Hurlbert \cite{glennsurvey}, and a survey of open problems in graph pebbling may be found at \cite{glennpage}.

Consider also the following {variant} of pebbling called
cover pebbling, first discussed by Crull et al (\cite{crull}):
The {\it cover pebbling number}
$\lambda(G)$ is defined as the minimum number of pebbles required such that it is possible, given any initial configuration of at least $\lambda(G)$ pebbles on $G$, to make a series of pebbling moves that {\it simultaneously} reaches {\it each} vertex of $G$.  A configuration is said to be {\it cover solvable} if it is
possible to place a pebble on every vertex of $G$ starting with that
configuration. Various results on cover pebbling have been determined. For
instance, we now know (\cite{crull}) that $\l(K_n)=2n-1; \l(P_n)=2^n-1$; and that for trees $T_n$,
\begin{equation}
\l(T_n)=\max_{v\in V(T_n)}\sum_{u\in V(T_n)}2^{d(u,v)},\end{equation}
where $d(u,v)$ denotes the distance between vertices $u$ and $v$.  Likewise, it was shown in \cite{munyon} that $\l(Q^d)=3^d$ and in \cite{nate} that $\l(K_{r_1,\ldots,r_m})=4r_1+2r_2+\ldots+2r_m-3$, where $r_1\ge3$ and $r_1\ge r_2\ge \ldots \ge r_m$.  The above examples reveal that for these special classes of graphs at any rate, the cover pebbling number equals the ``stacking number", or, put another way, the worst possible distribution of pebbles consists of placing all the pebbles on a single vertex.  The intuition built by computing the value of the cover pebbling number for the families $K_n$, $P_n$, and $T_n$ in \cite{crull} led to Open Question No.~10 in \cite{crull},  and which was proved by Vuong and Wyckoff \cite{ian}  and later, independently, by Sj\"ostrand \cite{jonas}:

\medskip

\noindent STACKING THEOREM:  {\it For any connected graph $G$,
$$\l(G)=\max_{v\in V(G)}\sum_{u\in V(G)}2^{d(u,v)},$$}

\medskip
\noindent thereby proving that (1) holds for all graphs.

The Sierpi\'nski graph will now be revealed to be one for which the use of the Stacking Theorem does not reduce the computation of the cover pebbling number to a trivial exercise.  We first prove that the diameter of $S_n$ is $2^{n-1}$, and, using this fact, that the worst vertex on which to stack pebbles is a corner vertex $a$ of degree two:

\begin{lemma}
$\diam(S_n)=2^{n-1}.$
\end{lemma}

 \medskip

\noindent {\bf Proof}  The fact that $\delta_n\ge2^{n-1}$ is obvious.  We use induction for the reverse inequality.  The result is clearly true for $n=1$.  Assume it to be true for $n$.  Let $x,y$ be any two points in (without loss of generality) $S_{n+1,T}$ and $S_{n+1,L}$ respectively.  Any path between $x$ and $y$ must go through $S_{n+1,L,T}$ which we will denote for brevity by $z$.  We thus have 
\begin{eqnarray*}
d(x,y)&=&d(x,z)+d(z,y)\\
&\le&2^{n-1}+2^{n-1}=2^n,
\end{eqnarray*}
as required.\hfill\qed

\begin{lemma} $ST(a)\ge ST(v)$ for each $v\in S_n$, where $ST(v)=\sum_{u\in S_n}2^{d(u,v)}$.
\end{lemma}

\medskip

\noindent{\bf Proof.} We proceed by induction.  The result is easy to verify for $S_1$ and $S_2$.  Assume that it is true for $S_n$. In $S_{n+1}$, denote the vertices $S_{n+1,T,T}, S_{n+1,L,T},$ and $S_{n+1,R,T}$ by $a$, $b$ and $c$ respectively.  Let $d$ and $e$ be arbitrary vertices in $S_{n+1,T}$ and $S_{n+1,L}\cup S_{n+1,R}$ respectively.  By Lemma 8, we have
\[d(a,e)={2^{n-1}}+\min\{d(b,e),d(c,e)\},\]
and
\[d(d,e)=\min\{d(d,b)+d(b,e), d(d,c)+d(c,e)\}.\]
Since, however, $\min\{d(b,e),d(c,e)\}=d(\alpha,e)$, where $\alpha=b$ or $\alpha=c$, it follows that
\[d(d,e)\le d(d,\alpha)+d(\alpha,e)\le{2^{n-1}}+d(\alpha,e)=d(a,e).\]
We thus have
\begin{eqnarray*}
ST(d)&=&\sum_{u\in S_{n+1}}2^{d(d,u)}\\
&=&\sum_{u\in S_{n+1,T}}2^{d(d,u)}+\sum_{u\in [S_{n+1,L}\cup S_{n+1,R}]\setminus\{b,c\}}2^{d(d,u)}\\
&\le&\sum_{u\in S_{n+1,T}}2^{d(d,u)}+\sum_{u\in [S_{n+1,L}\cup S_{n+1,R}]\setminus\{b,c\}}2^{d(a,u)}\\
&\le&\sum_{u\in S_{n+1,T}}2^{d(a,u)}+\sum_{u\in [S_{n+1,L}\cup S_{n+1,R}]\setminus\{b,c\}}2^{d(a,u)}\\
&=&ST(a),
\end{eqnarray*}
where the next to last line above follows due to the induction hypothesis.\hfill\qed

\begin{thm} The cover pebbling number $\lambda(S_n)$ of the Sierpi\'nski graph satisfies the recursion
\[\l(S_{n+1})=\lr1+2^{2^{n-1}+1}\rr\l(S_n)-\lr2^{2^n}+2^{2^{n-1}+1}\rr.\]
\end{thm}

\medskip

\noindent{\bf Proof} Notice that, with $\beta(i)$ denoting (in $S_{n+1}$) the number of points at distance $i$ from the corner vertex $a$, the following two conditions hold: 
\begin{itemize}
\item For any $j\in\{1,2,\ldots,2^{n-1}-1\}$, $\beta(j+2^{n-1})=2\beta(j)$;
\item $\beta(2^n)=2\beta(2^{n-1})-1$.  
\end{itemize}
Thus,
\begin{eqnarray*}
\l(S_{n+1})&=&ST(a)\\
&=&\sum_{i=0}^{2^{n-1}}\beta(i)2^i+2\sum_{i=1}^{2^{n-1}}\beta(i)2^{i+2^{n-1}}-2^{2^n}\\
&=&\l(S_n)+2\cdot2^{2^{n-1}}(\l(S_n)-1)-2^{2^n}\\
&=&\lr1+2^{2^{n-1}+1}\rr\l(S_n)-\lr2^{2^n}+2^{2^{n-1}+1}\rr,
\end{eqnarray*}
as asserted.\hfill\qed

\section{Open Problems} Here are some open problems for readers of this paper to consider:
\begin{itemize}
\item What is the edge chromatic number (chromatic index) and total chromatic number of $S_n$?  By Vizing's Theorem, the former is either 4 or 5, and if Behzad's total chromatic conjecture is true, then the latter is either 5 or 6.
\item What is the pebbling number of $S_n$?  Various bounds as in \cite{chan} may be used to estimate this quantity, but we consider the determination of $\pi(S_n)$ to be quite hard.
\item What baseline structural properties similar to the ones we have studied in this paper can be established for the Sierpi\'nski graphs of Klav\v zar and Milutinovi\'c \cite{klavzar}?
\item In a similar vein, what can be said of the domination number, cycle structure, etc.~of Sierpi\'nski-like graphs generated by considering Pascal's triangle mod $p;p\ge3$? (See \cite{shannon} for details on this structure and recall that the Siepi\'nski gasket graph is related to Pascal's triangle mod 2.)
\end{itemize}
\section{Acknowledgment}  Teguia was a graduate student at ETSU when this research was conducted --  but not as part of his M.S. thesis, which was in operator theory.  Godbole's research was  supported by NSF Grant DMS-0139291.  This paper has benefited greatly from the suggestions for improvement made by the anonymous referee.

\end{document}